\documentstyle[twoside]{article}
\include{graphicx}
\title{\bf Asymptotics of a ${}_3F_2$ hypergeometric function with four large parameters}
\author{\sc R. B. Paris\footnote{E-mail address:\ \ {\tt r.paris@abertay.ac.uk}}\\
\\
{\em Division of Computing and Mathematics,}\\
{\em Abertay University, Dundee DD1 1HG, UK}\\
}

\begin{document}
\newcommand{\bee}{\begin{equation}}
\newcommand{\ee}{\end{equation}}
\def\f#1#2{\mbox{${\textstyle \frac{#1}{#2}}$}}
\def\dfrac#1#2{\displaystyle{\frac{#1}{#2}}}
\newcommand{\fr}{\frac{1}{2}}
\newcommand{\fs}{\f{1}{2}}
\newcommand{\g}{\Gamma}
\newcommand{\br}{\biggr}
\newcommand{\bl}{\biggl}
\newcommand{\ra}{\rightarrow}
\renewcommand{\topfraction}{0.9}
\renewcommand{\bottomfraction}{0.9}
\renewcommand{\textfraction}{0.05}
\newcommand{\mcol}{\multicolumn}
\date{}
\maketitle
\pagestyle{myheadings}
\markboth{\hfill {\it R.B. Paris} \hfill}
{\hfill {\it Hypergeometric function with large parameters} \hfill}
\begin{abstract} 
We consider the asymptotic behaviour of the generalised hypergeometric function
\[{}_3F_2\bl(\!\!\begin{array}{c} 1, \fs(1+t)k, \fs(1+t)k+\fs\\tk+1, k+1\end{array}\!\!; x\br),\qquad 0<x,t\leq 1\]
as the parameter $k\to+\infty$.   Numerical results illustrating the accuracy of the resulting expansion are given. 
\vspace{0.4cm}

\noindent {\bf MSC:} 33C05, 34E05, 41A60
\vspace{0.3cm}

\noindent {\bf Keywords:} Hypergeometric function, asymptotic expansion, large parameters\\
\end{abstract}

\vspace{0.2cm}

\noindent $\,$\hrulefill $\,$

\vspace{0.2cm}

\begin{center}
{\bf 1. \  Introduction}
\end{center}
\setcounter{section}{1}
\setcounter{equation}{0}
\renewcommand{\theequation}{\arabic{section}.\arabic{equation}}
The following problem arising in two-variable moment theory \cite{Ger} is the determination of the asymptotic behaviour of the generalised hypergeometric function 
\bee\label{e11}
S(x;t):={}_3F_2\bl(\!\!\begin{array}{c} 1, \fs(1+t)k, \fs(1+t)k+\fs\\tk+1, k+1\end{array}\!\!; x\br), \qquad 0<x,t\leq 1
\ee
as the parameter $k\to+\infty$. The parametric excess\footnote{The parametric excess equals the difference between the sums of the denominator and numerator parameters.} of this function equals $\fs$ so that $S(x;t)$ converges as $x\to 1$. We also have the simple evaluation
\[S(1;0)={}_2F_1\bl(\!\!\begin{array}{c} \fs k, \fs k+\fs\\k+1\end{array}\!\!;1\br)=
\frac{\g(k+1)\g(\fs)}{\g(\fs k+\fs)\g(\fs k+1)}=2^k\]
by the Gauss summation formula.

An integral representation for $S(x;t)$ involving the modified Bessel function $K_\nu(z)$ can be obtained from \cite[p.~156]{S} in the form
\[S(x;t)=\frac{2}{\g(\nu)}\int_0^\infty u^{(\nu-1)/2} K_{\nu-1}(2\sqrt{u})\,{}_1F_2\bl(\!\!\begin{array}{c}\fs(1+t)k+\fs\\1+tk, k+1\end{array}\!\!;xu\br)du,\]
where $\nu=\fs(1+t)k$.

In this note we investigate the large-$k$ behaviour of $S(x;t)$. We first employ a contour integral representation for $S(x;t)$ that involves the Gauss hypergeometric function with two large parameters. The known asymptotics of this latter function then lead to an expansion for $S(x;t)$ also in terms of Gauss hypergeometric functions of a different form. The asymptotic expansion of these functions has recently been considered in \cite{P2}. In the second approach we make use of the confluence principle as discussed in \cite[pp.~56-57]{L}.
Numerical results are presented to demonstrate the accuracy of the expansion obtained.
\vspace{0.6cm}

\begin{center}
{\bf 2. \ An expansion for $S(x;t)$}
\end{center}
\setcounter{section}{2}
\setcounter{equation}{0}
\renewcommand{\theequation}{\arabic{section}.\arabic{equation}}
We employ the integral representation 
\bee\label{e21}
S(x;t)=\frac{\g(tk+1)\g(bk)}{\g(ak)}\,\frac{1}{2\pi i}\int_0^{(1+)} u^{ak-1}(u-1)^{-bk}\,{}_2F_1\bl(\!\!\begin{array}{c}1, ak+\fs\\k+1\end{array}\!\!;xu\br)du,
\ee
where we define
\bee\label{e21a}
a=\fs(1+t),\quad b=\fs(1-t),\quad c=a(1-a)=\f{1}{4}(1-t^2) \qquad (0<t<1).
\ee
The validity of this representation requires $b>0$, which implies that $t$ must satisfy $0<t<1$.
The integration path is a closed loop that starts at the origin, encircles $t=1$ in the positive sense (excluding the point $t=1/x$) and returns to the origin. The integral (\ref{e21}) can be established from the representation over $[0,1]$ for ${}_3F_2(x)$ given in \cite[(16.5.2)]{DLMF}
extended into a contour integral by means of \cite[(5.12.10)]{DLMF}.

In (\ref{e21}) the integration path can be arranged so that $|xu|<1$ everywhere on the loop.
For $k\to+\infty$, we can then employ the expansion of the Gauss hypergeometric function appearing in (\ref{e21}). This has two large parameters with $a<1$; from \cite[Section 3.1]{P1}, this function is associated with the Laplace integral
\[\int_0^1 f(\tau) e^{k\psi(\tau)} d\tau,\qquad \psi(\tau):=a\log\,\tau+(1-a)\log (1-\tau),\quad f(\tau):=\frac{\tau^{-1/2}(1-\tau)^{-1/2}}{1-z\tau}.\]
Expansion of this integral by Laplace's method about the saddle point $\tau=a$ yields the result \cite[Eq.~(3.6)]{P1}
\bee\label{e22}
{}_2F_1\bl(\!\!\begin{array}{c}1, ak+\fs\\k+1\end{array}\!\!;z\br)\sim \frac{\Xi(a,k)}{1-az}\bl\{1+\frac{c_2(a,z)}{k}+\frac{3c_4(a,z)}{4k^2}+\cdots\br\} \qquad (a<1),
\ee
where 
\[\Xi(a,k)=\frac{\g(k+1)}{\g(ak+\fs)\g((1-a)k+\fs)}\,\bl(\frac{k}{2\pi}\br)^{\!\!-1/2} a^{ak}(1-a)^{(1-a)k}.\]

The coefficients $c_2(a,z)$ and $c_4(a,z)$ are defined by
\[c_2(a,z)=-\frac{1}{\psi''}\{{F}_2-\Psi_3{F}_1+\f{5}{12}\Psi_3^2-\f{1}{4}\Psi_4\},\]
\[c_4(a,z)=\frac{1}{(\psi'')^2}\{\f{1}{6}{ F}_4-\f{5}{9}\Psi_3{F}_3+\f{5}{12}(\f{7}{3}\Psi_3^2-\Psi_4){ F}_2-\f{35}{36}(\Psi_3^3-\Psi_3\Psi_4+\f{6}{35}\Psi_5){F}_1 \]
\bee\label{e23}
+\f{35}{36}(\f{11}{24}\Psi_3^4-\f{3}{4}(\Psi_3^2-\f{1}{6}\Psi_4)\Psi_4+\f{1}{5}\Psi_3\Psi_5-\f{1}{35}\Psi_6)\},
\ee
where, for brevity, we have defined
\bee\label{e23a}
\Psi_k:=\frac{\psi^{(k)}(a)}{\psi''(a)}\ \ \ (k\geq 3),\qquad {F}_k:=\frac{{ f}^{(k)}(a)}{{f}(a)}\ \ \ (k\geq 1);
\ee
see, for example, \cite[p.~119]{D}, \cite[p.~127]{O} or \cite[p.~13]{PBook}. 
Substitution of the above forms of $\psi(\tau)$ and $f(\tau)$ into (\ref{e23}) and (\ref{e23a}) yields after some laborious algebra (carried out with the aid of {\it Mathematica}) the coefficient values given by
\[c_2(a,z)=\frac{2cz^2}{(1-az)^2}+\frac{(1-2a)z}{1-az}-\frac{1+2c}{12c},\]
\[c_4(a,z)=\frac{4c^2z^4}{(1-az)^4}+\frac{14(1-2a)cz^3}{3(1-az)^3}+\frac{(3-20c)z^2}{3(1-az)^2}-\frac{2(1-2a)z}{3(1-az)}-\frac{(1+2c)^2}{864c^2}\hspace{2cm}\]
\bee\label{e26}
\hspace{8cm}-\frac{(1+2c)}{36c}\,c_2(a,z)~,
\ee
where the quantity $c$ is specified in (\ref{e21a}).

Application of Stirling's formula for the gamma function shows that
\bee\label{e24}
\Xi(a,k)=1+\frac{1+2c}{24ck}+\frac{(1+2c)^2}{1152c^2k^2}+O(k^{-3}).
\ee
as $k\to+\infty$. Then with $z=xu$ in (\ref{e26}) and use of (\ref{e24}) we obtain
\[S(x;t)\sim\frac{\g(tk+1)\g(bk)}{\g(ak)}\,\frac{\Xi(a,k)}{2\pi i}\int_0^{(1+)}\frac{u^{ak-1}(u-1)^{-bk}}{1-axu}\bl\{1+\frac{c_2(a,xu)}{{2k}}+\frac{3c_4(a,xu)}{4k^2}+\cdots\br\}du\]
\bee\label{e25}
=\frac{\g(tk+1)\g(bk)}{\g(ak)}\,\frac{1}{2\pi i}\int_0^{(1+)} \frac{u^{ak-1}(u-1)^{-bk}}{1-axu}\bl\{\Xi(a,k)+\frac{c_2(a,xu)}{2k}+\frac{3c_4'(a,xu)}{4k^2}+\cdots\br\}du,
\ee
where
\[c_4'(a,xu)=c_4(a,xu)+\frac{(1+2c)}{36c}\,c_2(a,xu).\]

We now employ the result \cite[(15.6.2)]{DLMF}
\[\frac{1}{2\pi i}\int_0^{(1+)} \frac{u^{\beta-1}(u-1)^{\gamma-\beta-1}}{(1-uz)^\alpha}\,du=
\frac{\g(\beta)}{\g(\gamma)\g(1+\beta-\gamma)}\,{}_2F_1(\alpha,\beta;\gamma;z)\]
when $\Re (\beta)>0$ and $\gamma-\beta\neq 1, 2, 3, \ldots\,$, 
to find that, for $m=0, 1, 2, \ldots\,$,
\[\frac{\g(tk+1)\g(bk)}{\g(ak)}\,\frac{1}{2\pi i}\int_0^{(1+)} \frac{u^{ak+m-1}(u-1)^{-bk}}{(1-ax u)^{m+1}}\,du=A_m\,{\cal F}_m,\]
where 
\bee\label{e25c}
A_m:=\frac{(ak)_m}{(tk+1)_m},\qquad{\cal F}_m:={}_2F_1\bl(\begin{array}{c}m+1, ak+m\\tk+m+1\end{array}\!\!;ax\br).
\ee
The terms involving $(1-axu)^{-1}$ in (\ref{e25}) combine to yield
\[\bl\{\Xi(a,k)-\frac{1+2c}{24ck}-\frac{(1+2c)^2}{1152c^2k^2}\br\} {\cal F}_0={\cal F}_0 \{1+O(k^{-3})\}\qquad(k\to+\infty)\]
by (\ref{e24}). 

Evaluation of the remaining terms then produces the following expansion
\[S(x;t)\sim {\cal F}_0+\frac{1}{k}\bl\{\frac{1}{2}(1-2a)xA_1{\cal F}_1+a(1-a)x^2A_2{\cal F}_2\br\}\]
\bee\label{e25d}
+\frac{1}{k^2}\bl\{\frac{1}{2}(2a-1)xA_1{\cal F}_1+\frac{1}{4}(3-20c)x^2A_2{\cal F}_2+\frac{7}{2}(1-2a)cx^3A_3{\cal F}_3+3c^2x^4A_4{\cal F}_4\br\}+\ldots\,
\ee
as $k\to+\infty$.
\vspace{0.6cm}

\begin{center}
{\bf 3. \ An alternative approach}
\end{center}
\setcounter{section}{3}
\setcounter{equation}{0}
\renewcommand{\theequation}{\arabic{section}.\arabic{equation}}
We give an alternative derivation of the expansion (\ref{e25d}) based on the confluence principle described in 
\cite[pp.~56--57]{L}. This approach is valid for $0<t\leq 1$, so that $a\leq 1$. Proceeding formally, we have the series representation
\bee\label{e31}
S(x;t)=\sum_{r\geq 0}\frac{(1)_r (ak)_r(ak+\fs)_r}{(tk+1)_r(k+1)_r}\,\frac{x^r}{r!}=\sum_{r\geq 0}\frac{(1)_r (ak)_r}{(tk+1)_r}\,\frac{(ax)^r}{r!}\,P_r,
\ee
where
\[P_r=\frac{(ak+\fs)_r}{(k+1)_r}=a^{-r}\frac{(1+\frac{1}{2ak})(1+\frac{3}{2ak})\cdots (1+\frac{2r-1}{2ak})}{(1+\frac{1}{k})(1+\frac{2}{k})\cdots (1+\frac{r}{k})}.\]
It follows that, for $k\to+\infty$,
\begin{eqnarray*}
\log\,P_r&=&\frac{1}{2ak}\sum_{n=1}^r (2n-1)-\frac{1}{8a^2k^2}\sum_{n=1}^r(2n-1)^2-\frac{1}{k}\sum_{n=1}^rn+\frac{1}{2k^2}\sum_{n=1}^rn^2+O(k^{-3})\\
&=&\frac{r^2}{2ak}-\frac{r(4r^2-1)}{24a^2k^2}-\frac{r(r+1)}{2k}+\frac{r(r+1)}{12k^2}(2r+1)+O(k^{-3})\\
&=&\frac{r^2(1-a)-ar}{2ak}+\frac{1}{24a^2k^2}\{(2a^2+1)r+6a^2r^2-4(1-a^2)r^3\}+O(k^{-3}),
\end{eqnarray*}
which upon exponentiation then yields
\[P_r=1-\frac{ar-r^2(1-a)}{2ak}+\frac{1}{24a^2k^2}\{(2a^2+1)r+9a^2r^2-2(1-a)(5a+2)r^3\]
\[\hspace{7cm}+3(1-a)^2r^4\}+O(k^{-3}).\]
This expansion assumes that $r^2/k\ll 1$. But the summation index $r\in[0,\infty)$; if the terms in the sum (\ref{e31}) are negligible for $r\gg r_0$, then if $k$ is such that $r_0^2/k\ll 1$ we can formally neglect the tail $r\gg r_0$, substitute the expansion for $P_r$ into (\ref{e31}) and evaluate the resulting series term by term.

We now employ the result valid for convergent series
\[\sum_{r\geq 0}b_r r^m \chi^r=\Theta^m \sum_{r\geq 0}b_r \chi^r,\qquad \Theta:=\chi \frac{d}{d\chi}\qquad(m=0, 1, 2, \ldots),\]
where we identify the coefficients $b_r$ with $(ak)_r/(tk+1)_r$. Then substitution of the above expansion for $P_r$ into (\ref{e31}) leads to
\[S(x;t)\sim {\cal F}_0+\frac{1}{2ak}\{-a\Theta {\cal F}_0+(1-a)\Theta^2 {\cal F}_0\}\]
\[+\frac{1}{24a^2k^2}\{(2a^2+1)\Theta{\cal F}_0+9a^2\Theta^2{\cal F}_0-2(1-a)(2+5a)\Theta^3{\cal F}_0+3(1-a)^2\Theta^4{\cal F}_0\}+\ldots\, ,\]
where
\[{\cal F}_0=\sum_{r\geq 0}\frac{(1)_r(ak)_r}{(tk+1)_r}\,\frac{\chi^r}{r!}={}_2F_1\bl(\begin{array}{c}1, ak\\tk+1\end{array}\!\!;\chi\br),\quad \chi:=ax.\]
Since
\[\Theta{\cal F}_0=\chi{\cal F}_0',\quad \Theta^2{\cal F}_0=\chi^2{\cal F}_0''+\chi{\cal F}_0',\quad\Theta^3{\cal F}_0=\chi^3{\cal F}_0'''+3\chi^2{\cal F}_0''+\chi{\cal F}_0',\]
\[\Theta^4{\cal F}_0=\chi^4{\cal F}_0^{(iv)}+6\chi^3{\cal F}_0'''+7\chi^2{\cal F}_0''+\chi{\cal F}_0',\]
we find
\[S(x;t)\sim{\cal F}_0+\frac{1}{2ak}\{(1-2a)\chi{\cal F}_0'+(1-a)\chi^2{\cal F}_0''\}\]
\[+\frac{1}{24a^2k^2}\{12a(2a-1)\chi{\cal F}_0'+3(3-20c)\chi^2{\cal F}_0''+14(1-a)(1-2a)\chi^3{\cal F}_0'''+3(1-a)^2\chi^4{\cal F}_0^{(iv)}\}+\ldots\ .\]

From \cite[(15.5.2)]{DLMF}, the derivatives of ${\cal F}_0$ are given by
\[{\cal F}_0^{(m)}=\frac{d^m}{d\chi^m}\, {}_2F_1\bl(\begin{array}{c}1, ak\\tk+1\end{array}\!\!;\chi\br)=m! A_m{\cal F}_m,\]
where $A_m$ and ${\cal F}_m$ are defined in (\ref{e25c}). This enables the expansion of $S(x;t)$ to be written in the form
\[
S(x;t)\sim{\cal F}_0+\frac{1}{k}\bl\{\frac{1}{2}(1-2a)xA_1{\cal F}_1+a(1-a)x^2A_2{\cal F}_2\br\}\]
\bee\label{e32}
+\frac{1}{k^2}\bl\{\frac{1}{2}(2a-1)xA_1{\cal F}_1+\frac{1}{4}(3-20c)x^2A_2{\cal F}_2+\frac{7}{2}(1-2a)cx^3A_3{\cal F}_3+3c^2x^4A_4{\cal F}_4\br\}+\ldots\,,
\ee
where we recall that $c=a(1-a)$. This expansion agrees with that in (\ref{e25d}) obtained from the contour integral
approach. We remark that its derivation has allowed us to consider $a\leq 1$ ($0<t\leq 1$) whereas the integral (\ref{e21}) requires $b>0$ ($0<t<1$).

\vspace{0.6cm}

\begin{center}
{\bf 4. \ Numerical results}
\end{center}
\setcounter{section}{4}
\setcounter{equation}{0}
\renewcommand{\theequation}{\arabic{section}.\arabic{equation}}
In this section we expand the quantities $A_m$ appearing in (\ref{e25d}) in inverse powers of the large parameter $k$ to obtain a modified expansion for $S(x;t)$.

For $k\to+\infty$ (with $t$ bounded away from zero) we have
\[A_1=\frac{a}{t}\bl\{1-\frac{1}{kt}+O(k^{-2})\br\},\qquad A_2=\frac{a^2}{t^2}\bl\{1-\frac{\alpha}{4k}+O(k^{-2})\br\}\]
with
\[\alpha=4\bl(\frac{3}{t}-\frac{1}{a}\br)=\frac{4(1+a)}{at},\]
and generally $A_m=(a/t)^m\{1+O(k^{-1})\}$. Then, upon expanding the quantities $A_m$ ($1\leq m\leq 4$) and noting that $2a-1=t$, we obtain from (\ref{e25d}) or (\ref{e32}) the expansion in the following modified form
\[S(x;t)\sim {\cal F}_0-\frac{1}{k}\bl\{\frac{1}{2}tX{\cal F}_1-cX^2{\cal F}_2\br\}\hspace{5cm}\]
\bee\label{e41}
+\frac{1}{k^2}\bl\{aX{\cal F}_1+\frac{1}{4}[3-(20+\alpha)c\,]X^2{\cal F}_2-\frac{7}{2}ctX^3{\cal F}_3+3c^2X^4{\cal F}_4\br\}+\ldots \, ,
\ee
where we have put $X:=ax/t$. This expansion holds for $k\to+\infty$, provided $t$ is bounded away from zero such that $kt\to+\infty$.
In the case $t=1$ (so that $a=1$, $c=0$), (\ref{e41}) reduces to
\bee\label{e42}
S(x;1)\sim {\cal F}_0-\frac{X{\cal F}_1}{2k}
+\frac{1}{k^2}\bl\{X{\cal F}_1+\frac{3}{4}X^2{\cal F}_2\br\}+\ldots \qquad(k\to+\infty).
\ee

The leading term of the above expansions is given by
\[{\cal F}_0={}_2F_1\bl(\begin{array}{c}1, ak\\tk+1\end{array}\!\!;ax\br).\]
Since the parametric excess associated with this function equals $-bk<0$ the series converges if $ax<1$. Thus ${\cal F}_0$ converges for $x\leq 1$ if $t<1$ and for $x<1$ if $t=1$. In Table 1 we present values of the absolute relative error in the computation of $S(x;t)$ using the expansions (\ref{e41}) and (\ref{e42}) when truncated at the term $O(k^{-M})$, $0\leq M\leq 2$.
\begin{table}[th]
\caption{\footnotesize{Values of the absolute relative error in the computation of $S(x;t)$ from (\ref{e41}) and (\ref{e42}) as a function of the truncation index $M$ and different values of $k$, $x$  and $t$.}} \label{t1}
\begin{center}
\begin{tabular}{|l||c|c|c|c|}
\hline
%&&&&\\[-0.3cm]
\mcol{1}{|c}{} & \mcol{2}{c|}{$k=100$} & \mcol{2}{c|}{$k=200$}\\
\mcol{1}{|c||}{$M$} & \mcol{1}{c|}{$x=0.50, t=0.75$} & \mcol{1}{c|}{$x=0.50, t=1$} & \mcol{1}{c|}{$x=0.75, t=0.50$}
& \mcol{1}{c|}{$x=0.50, t=0.50$}\\
\hline
&&&&\\[-0.25cm]
0 & $5.723\times 10^{-3}$ & $9.481\times 10^{-3}$ & $1.357\times 10^{-1}$ & $9.638\times 10^{-4}$ \\
1 & $9.925\times 10^{-5}$ & $3.288\times 10^{-4}$ & $6.171\times 10^{-4}$ & $1.286\times 10^{-4}$ \\
2 & $1.223\times 10^{-6}$ & $1.315\times 10^{-5}$ & $2.380\times 10^{-4}$ & $4.238\times 10^{-6}$ \\
[.1cm]\hline
\mcol{1}{|c}{} & \mcol{2}{c|}{$k=200$} & \mcol{2}{c|}{$k=300$}\\
\mcol{1}{|c||}{$M$} & \mcol{1}{c|}{$x=0.50, t=0.75$} & \mcol{1}{c|}{$x=0.50, t=1$} & \mcol{1}{c|}{$x=0.75, t=0.50$}
& \mcol{1}{c|}{$x=0.50, t=0.50$}\\
\hline
&&&&\\[-0.25cm]
0 & $2.919\times 10^{-3}$ & $4.866\times 10^{-3}$ & $1.073\times 10^{-1}$ & $7.293\times 10^{-4}$ \\
1 & $2.458\times 10^{-5}$ & $8.476\times 10^{-5}$ & $1.150\times 10^{-3}$ & $5.973\times 10^{-5}$ \\
2 & $1.897\times 10^{-7}$ & $1.710\times 10^{-6}$ & $6.018\times 10^{-5}$ & $1.285\times 10^{-6}$ \\
[.1cm]\hline
\end{tabular}
\end{center}
\end{table}

If we write $\epsilon=a/t$, $\lambda=tk$ and $\chi=ax$, we have when $0<t<1$
\bee\label{e43}
{\cal F}_0={}_2F_1\bl(\begin{array}{c}1,  \epsilon\lambda\\1+\lambda\end{array}\!\!;\chi\br),\qquad \epsilon>1.
\ee
The expansion of the Gauss hypergeometric function ${}_2F_1(1, \alpha+\epsilon\lambda;\beta+\lambda;\chi)$ for $\epsilon>0$ and $\lambda\to+\infty$ has been considered in \cite{P1}, and more recently when $\epsilon>1$ in \cite{P2}. The analysis of this function when $\epsilon>1$ is complicated by the presence in its Laplace-type integral representation (where a multiplicative factor is omitted) 
\[\frac{1}{2\pi i}\int_0^{(1+)} \frac{\tau^{\alpha-1}(\tau-1)^{\beta-\alpha-1}}{1-\chi\tau}\,e^{-\lambda\phi(\tau)}d\tau,\qquad \phi(\tau)=(\epsilon-1) \log (\tau-1)-\epsilon\log\,\tau\]
of a saddle point at $\tau=\epsilon$ and a simple pole at $\tau=1/\chi$. These points coalesce when $\epsilon\chi=1$; that is, since $\epsilon\chi=a^2x/t$, when 
\[x=x^*=\frac{4t}{(1+t)^2}.\]
From \cite[Eq.~(3.13)]{P1}, we have as $\lambda\to+\infty$ (when $\epsilon\chi<1$)
\bee\label{e44}
{\cal F}_0\sim
%{}_2F_1\bl(\begin{array}{c}1, \epsilon\lambda\\1+\lambda\end{array}\!\!;\chi\br)\sim
\frac{G(\lambda)}{\sqrt{2\pi}}\,\bl(\frac{\epsilon}{\epsilon-1}\br)^{\!\!-1/2}\,\frac{e^{-\lambda\phi(\epsilon)}}{1-\epsilon\chi}\sum_{k=0}^\infty \frac{c_{2k} \g(k+\fs)}{\lambda^{k+\frac{1}{2}}\g(\fs)}\qquad(\epsilon>1,\  x<x^*),
\ee
where
\[G(\lambda):=\frac{\g(1+\lambda) \g((\epsilon-1)\lambda)}{\g(\epsilon\lambda)}\sim(2\pi\lambda)^{1/2}\bl(\frac{\epsilon}{\epsilon-1}\br)^{\!\!1/2} e^{\lambda\phi(\epsilon)}\qquad(\lambda\to+\infty)\]
by application of Stirling's formula.
The coefficient $c_0=1$ with $c_2$, $c_4$ given in \cite[Eqs.~(3.4), (3.5)]{P1}. This shows that the leading large-$\lambda$ behaviour of ${\cal F}_0$ when $x<x^*$ is given by $(1-\epsilon\chi)^{-1}$. 

The expansion (\ref{e44}) breaks down as $x\to x^*$. From \cite[Theorem 1]{P2} we have the uniform expansion (corresponding to $\alpha=0$, $\beta=1$)
\bee\label{e45}
{\cal F}_0\sim
%{}_2F_1\bl(\begin{array}{c}1, 1+\epsilon\lambda\\1+\lambda\end{array}\!\!;\chi\br)\sim
\frac{G(\lambda)}{2}\bl\{\chi e^{-\lambda\phi(1/\chi)}\,\mbox{erfc}\,(\pm\lambda^\frac{1}{2}p)+\frac{e^{-\lambda\phi(\epsilon)}}{\pi\chi}\sum_{k=0}^\infty d_{2k} \frac{\g(k+\fs)}{\lambda^{k+\frac{1}{2}}}\br\},
\ee
where erfc denotes the complementary error function and
$p=(\phi(\epsilon)-\phi(1/\chi))^{1/2}\geq 0$. 
A closed-form expression for the coefficients $d_{2k}$  is given in \cite[Eq.~(2.17)]{P2} and 
\bee\label{e46}
d_0=\bl(\frac{2(\epsilon-1)}{\epsilon}\br)^{\!\!1/2}\!\!\frac{\chi}{1-\epsilon\chi}\mp\frac{\chi}{p}~.
\ee
The upper signs in (\ref{e45}) and (\ref{e46}) apply when $\epsilon\chi<1$ ($x<x^*$) and the lower signs when $\epsilon\chi>1$ ($x^*<x<1$).
These coefficients possess a removable singularity when the saddle and pole coincide, corresponding to $\epsilon\chi=1$ ($x=x^*$) and $p=0$; see \cite[Section 3]{P2} for the expansion at coalescence and an explicit representation of the first three coefficients. In Table 2 we show values of the absolute relative error in the computation of ${\cal F}_0$ from (\ref{e45}) when $t=1/3$, $k=150$ (so that $x^*=0.75$), corresponding to $\epsilon=2$,
$\chi=2x/3$ and $\lambda=50$, as a function of the truncation index $M$. The coefficients $d_{2k}$ for $k=1, 2$ are taken from \cite{P2}. 
\begin{table}[th]
\caption{\footnotesize{Values of the absolute relative error in the computation of ${\cal F}_0$ from (\ref{e45}) as a function of the truncation index $M$ when $t=1/3$ and $k=150$.}} \label{t2}
\begin{center}
\begin{tabular}{|l||l|l|l|l|l|}
\hline
&&&&&\\[-0.3cm]
\mcol{1}{|c||}{$M$} & \mcol{1}{c|}{$x=0.45$} & \mcol{1}{c|}{$x=0.72$} & \mcol{1}{c|}{$x=0.78$}
& \mcol{1}{c|}{$x=0.90$} & \mcol{1}{c|}{$x=1$}\\
\hline
&&&&&\\[-0.25cm]
0 & $6.025\times 10^{-5}$ & $1.353\times 10^{-6}$ & $7.122\times 10^{-7}$ & $3.455\times 10^{-7}$ &  $1.270\times 10^{-8}$  \\
1 & $6.803\times 10^{-7}$ & $1.664\times 10^{-8}$ & $8.762\times 10^{-9}$ & $4.168\times 10^{-9}$ & $1.474\times 10^{-10}$  \\
2 & $2.403\times 10^{-9}$ & $4.600\times 10^{-11}$ & $2.422\times 10^{-11}$ & $1.224\times 10^{-11}$& $4.835\times 10^{-13}$  \\
[.1cm]\hline
\end{tabular}
\end{center}
\end{table}

The expansion of the functions ${\cal F}_m$ with $m\geq 1$ can be obtained from that of ${}_2F_1(1,\alpha+\epsilon\lambda;\beta+\lambda;\chi)$, with $\alpha=m$, $\beta=m+1$, by means of a recurrence relation; see \cite[Section 4]{P2} for details.

\end{document}